\documentclass[11pt]{article}
\usepackage[active]{srcltx}
\usepackage{amssymb,amsmath}
\usepackage{color}
\usepackage[unicode]{hyperref}
\usepackage{txfonts}
\usepackage{graphicx}

\hypersetup{
    colorlinks = true,%
    citecolor = [rgb]{0.0,0.0,0.9},
    filecolor=black,%
    linkcolor = [rgb]{0.65,0.0,0.0},%
    anchorcolor = red,
    pagecolor = red,
    urlcolor= [rgb]{0.65,0.0,0.0}
}

\newtheorem{example}{ Example}[section]
\newtheorem{proposition}{Proposition}[section]
\newtheorem{theorem}{Theorem}[section]

\newtheorem{definition}{Definition}[section]
\newtheorem{remark}{Remark}[section]

\numberwithin{equation}{section}

\def\R{\mathbb{R}}

\def\psn{\par\medskip\noindent}

\def\qed{\hfill$\square$\psn\medskip}
\def\x{\overline{x}}

\begin{document}

\title{The Engulfing Property from a Convex Analysis Viewpoint}

\author{
 A. Calogero\thanks{
Dipartimento di Matematica e Applicazioni, Universit\`a degli Studi
di Milano-Bicocca, Via Cozzi 53, 20125 Milano, Italy ({\tt
andrea.calogero@unimib.it}, $^+$corresponding author: {\tt
rita.pini@unimib.it})},
\quad R. Pini$^*$$^+$
 }

\maketitle

\begin{abstract} In this note we provide a simple proof of some properties enjoyed by convex functions having the engulfing property. In particular, making use only of results peculiar to convex analysis, we prove that differentiability and strict convexity are conditions intrinsic to the engulfing property.
\end{abstract}



\section{Introduction}

The convex functions satisfying the so called \emph{engulfing property} have been studied in connection with the solution to the  Monge-Amp\`{e}re equation. Several conditions on such functions have been proposed in order to preserve the harmony between measure theory (and, in particular, the Monge-Amp\`{e}re measure related to a convex function), and the shape of the sections, with their induced geometry;
in this framework, we would like to mention  the celebrated $\mathcal{C}^{1,\beta}$-estimate due to Caffarelli [1-2],
and the exhaustive book by Guti\'errez [3]. This study involved many authors with different points of view; very interesting are the papers  by Guti\'errez and Huang [4], and by Forzani and Maldonado [5-6].
Let us devote a few lines to the mentioned $\mathcal{C}^{1,\beta}$-estimate: one of the proposed conditions on a convex function, with bounded sections, is the so called (DC)-doubling property of the related Monge-Amp\`{e}re measure (for details, see, for example, [6]); this condition plays a fundamental role  in the whole theory of the Monge-Amp\`{e}re equation, since it is equivalent to the engulfing property of the related function, and it implies that this function is strictly convex and in $ {\mathcal C}^{1,\beta}_{\rm loc}$. Due to this equivalence, the study of the engulfing property of a convex function is often moved to the study of the regularity of the related Monge-Amp\`{e}re measure.

The purpose of this note, differently from the literature, is to bring into focus that strict convexity and differentiability are properties intrinsic to the engulfing. This is done by means of purely convex analytic elementary
techniques, and without taking into account the properties of the related Monge-Amp\`{e}re measure.

The paper is organized as follows: in Section 2, we recall the notion of engulfing for a convex function, and introduce the notion of soft engulfing. We provide a fine monotonicity result for the subdifferential map of a convex function that enjoys the soft engulfing property. Furthermore, we prove that this property of the function entails continuous differentiability, as well as strict convexity. In Section 3, we prove the equivalence between the class of functions satisfying the engulfing property, and the class of functions satisfying the condition, apparently milder, of soft engulfing. Finally, in Section 4, some further directions of investigation are traced.

\section{Study of the Soft Engulfing Property of a Convex Function}

Given a convex function $\varphi:\R^n\to \R$, for every $x_0\in \R^n,$ $p\in \partial \varphi(x_0),$ and $s>0,$  we will denote by $S_{\varphi}(x_0,p,s)$ the section of $\varphi$ at $(x_0,\ p),$ with height $s$, defined as follows:
 \begin{equation}\label{def section R}
S_\varphi(x_0,p,s)=\{x\in \R^n:\ \varphi(x)- \varphi(x_0)-p\cdot(x-x_0)<s\}.
\end{equation}
In case $\varphi$ is differentiable at $x_0,$ we will denote  the section at $x_0,$ with height $s,$ by $S_\varphi(x_0,s),$ for short.
\begin{definition} We say that a convex function $\varphi$ satisfies
the \emph{engulfing property}
(shortly, $\varphi \in E(n,K),$ where $n$ denotes the dimension of the domain) if there exists $K>  1$ such that, for any $x\in\R^n,$ $p\in \partial \varphi(x),$ and
$t > 0$, if $y\in S_\varphi(x,p,t)$, then
$$
S_\varphi(x,p,t)\subset S_\varphi(y,q,Kt),
$$
for every $q\in \partial \varphi(y)$.
Likewise, we say that a convex function $\varphi$ satisfies
the \emph{soft engulfing property}
(shortly, $\varphi \in E^{\mathrm{soft}}(n,K)$)
if there exists $K>1$ such that, for any $x\in\R^n,$ $p\in \partial \varphi(x),$ and
$t > 0$, if $y\in S_\varphi(x,p,t)$, then
$$
x\in S_\varphi(y,q,Kt),
$$
for every $q\in \partial \varphi(y)$.
\end{definition}
Let us stress that, differently from the literature, in the previous definitions of the engulfing properties we do not require neither that the functions involved are differentiable, nor that their sections are bounded sets.
Trivially, one has that $E(n,K)\subset E^{\mathrm{soft}}(n,K);$ however, this \lq\lq soft\rq\rq\, condition is only apparently milder that the previous one, as we will see later on.

Let us recall some well-known notions and results. Given a multivalued map $T:\R^n\to {\mathcal P}(\R^n),$ we denote by
$\mathrm{dom}(T)$ and $\mathrm{gph}(T),$ respectively, the sets
$$
\mathrm{dom}(T):=\{x\in \R^n:\, T(x)\neq \emptyset\},
$$
and
$$
\mathrm{gph}(T):=\{(x,p)\in \R^n\times\R^n:\, p\in T(x)\}.
$$
Given a convex function $\varphi:\R^n\to\R$, the subdifferential $\partial \varphi(x_0)$ of $\varphi$ at $x_0$ is defined as the set
$$
\partial \varphi(x_0)=\{p\in \R^n:\, \varphi(x)\ge \varphi(x_0)+ p\cdot(x-x_0),\quad \forall x\in \R^n\}.
$$
The set $\partial \varphi(x_0)$ is nonempty, compact and convex, and, if $\varphi$ is differentiable at $x_0,$ then $\partial \varphi(x_0)=\{\nabla \varphi(x_0)\}.$
From a classical result in convex analysis, the convexity of $\varphi$ is fully characterized by the nonemptiness of $\partial\varphi(x),$ for every $x\in \R^n.$
The subdifferential map $\partial\varphi:\R^n\to {\mathcal P}(\R^n)$ is given by $x\mapsto \partial \varphi(x);$
its graph is a maximal monotone set, i.e., for every $(x,p),(y,q)\in \mathrm{gph}(\partial\varphi),$ we have that
$$
(p-q)\cdot(x-y)  \ge 0,
$$
and $\mathrm{gph}(\partial\varphi)$ cannot be extended to a monotone set without loosing the previous property
(see, for instance, Theorem 23.4 in [7]).

The celebrated and already mentioned result due to Caffarelli establishes that, given a Borel measure $\nu$ defined in $\R^n$, every strictly convex generalized solution (or Aleksandrov solution) $\varphi$ of the
Monge-Amp\`{e}re equation
\begin{equation}\label{Monge Ampere equation}
{\rm det}\, D^2\varphi=\nu
\end{equation}
 must be in the class ${\mathcal C}^{1,\beta}_{\rm loc},$ for some $\beta\in ]0,1[$. Let us recall briefly (see [3] for all the details) that the  Monge-Amp\`{e}re measure $\mu_\varphi$ associated to $\varphi$ is defined by
$$\mu_\varphi(E)=|\partial \varphi(E)|,$$
where the set
$\{E\subset \R^n:\, \partial\varphi(E)\;\mathrm{is}\;\mathrm{Lebesgue}\;\mathrm{measurable}\}$ is a $\sigma$-algebra containing the Borel sets, and $|\cdot|$ denotes  the Lebesgue measure; moreover,  $\varphi$ is a generalized solution of \eqref{Monge Ampere equation} if $\mu_\varphi=\nu$.

The study of the Monge-Amp\`{e}re measure $\mu_\varphi$ with the doubling property turns out to be intrinsically connected with the engulfing property of $\varphi$. Our point of view takes advantage of  convex analysis to investigate the properties of the functions $\varphi\in E(n,K)$ and, in particular, the behaviour of the sections  $S_{\varphi}(x_0,p,s),$ when we let $x_0$ and $s$ vary.

In the sequel, our aim will be to shed some light on the properties enjoyed by the functions in the class $E^{\mathrm{soft}}(n,K).$ The next proposition shows that the engulfing property of $\varphi$ is related to the monotonicity of the subdifferential map $\partial\varphi$. More precisely, it is known that, for every convex function $\varphi,$ the multivalued map $\partial\varphi$ is monotone; the following result highlights a finer behaviour of this map:

\begin{proposition}\label{prop:inequalities} Let $\varphi:\R^n\to \R$ be a function in $E^{\mathrm{soft}}(n,K).$ Then
\begin{equation}\label{caratterizzazione con p}
\begin{array}{ll}
\displaystyle\frac{K+1}{K}(\varphi(y)-\varphi(x)&-p\cdot(y-x))\le (p-q)\cdot(x-y)\\
&\le (K+1)\left(\varphi(y)-\varphi(x)-p\cdot(y-x)\right),\\
\end{array}
\end{equation}
for every $(x,p),(y,q)\in \mathrm{gph}(\partial \varphi).$
\end{proposition}
{\it Proof}\quad Take any $\epsilon >0;$ we have that
$$
\varphi(x)<\varphi(x)+\epsilon=\varphi(y)+q\cdot(x-y)+\varphi(x)-\varphi(y)-q\cdot(x-y)+\epsilon.
$$
Furthermore, $q\in \partial \varphi(y)$ implies that $\varphi(x)-\varphi(y)-q\cdot(x-y)\ge 0.$ Therefore,
$$
x\in S_\varphi(y,q,\varphi(x)-\varphi(y)-q\cdot(x-y)+\epsilon).
$$
By the soft engulfing property, $y\in S_\varphi(x,p,K(\varphi(x)-\varphi(y)-q\cdot(x-y)+\epsilon)),$ for every $p\in \partial \varphi(x),$ i.e.
$$
\varphi(y)<\varphi(x)+p\cdot(y-x)+K\varphi(x)-K\varphi(y)-Kq\cdot(x-y)+K\epsilon.
$$
Letting $\epsilon \searrow 0,$ we obtain
$$
(K+1)\varphi(y)\le (K+1)\varphi(x)+(p+Kq)\cdot(y-x),
$$
for every $(x,p),(y,q)\in \mathrm{gph}(\partial \varphi).$
By interchanging the roles of $x$ and $y,$ we get
$$
(K+1)\varphi(x)\le (K+1)\varphi(y)+(q+Kp)\cdot(x-y),
$$
for every $(x,p),(y,q)\in \mathrm{gph}(\partial \varphi).$
From the inequalities above we easily get
\begin{equation}\label{eq:1}
\frac{p+Kq}{K+1}\cdot(x-y)\le \varphi(x)-\varphi(y)\le \frac{q+Kp}{K+1}\cdot(x-y).
\end{equation}
From the first inequality in \eqref{eq:1}, it follows that
$$
\frac{1}{K+1}(p-q)\cdot(x-y)\le \varphi(x)-\varphi(y)-q\cdot(x-y).
$$
The second inequality in \eqref{eq:1} gives
$$
\varphi(x)-\varphi(y)-q\cdot(x-y)\le \frac{K}{K+1}(p-q)\cdot(x-y),
$$
thereby showing the assertion.
\qed

The property (\ref{caratterizzazione con p}) in Proposition \ref{caratterizzazione con p}, in fact, requires a sort of uniform control over the monotonicity of the subdifferential map in the whole $\R^n$. By means of (\ref{caratterizzazione con p}), it is an easy task to verify that the function $\varphi:\R\to\R$, defined by $\varphi(x)=e^x$, is not in $E^{\mathrm{soft}}(1,K),$ for any $K$; the reason of this exclusion does not depend on the asymmetry of this function, since one can easily prove that the function $x\mapsto e^{x^2}$ is not in  $E^{\mathrm{soft}}(1,K),$ too.

Furthermore, while every strictly convex polynomial enjoys the engulfing property (see, for instance,  [6]), not every strictly convex function with a \lq\lq polynomial behaviour\rq\rq\ does. The following example gives an idea to the non familiar reader:
\begin{example}\quad
\rm{Let us consider the function $\varphi:\R\to \R$ defined by
\begin{equation}\label{esempio}
\varphi(x)=
\left\{\begin{array}{ll}
x^4,& \;x\ge 0,\\
x^2,& \;x<0.\\
\end{array}
\right.
\end{equation}
We show that $\varphi\notin  E^{\mathrm{soft}}(1,K),$ for any $K>1.$ By contradiction, via Proposition \ref{prop:inequalities}, by taking  $x>0$ and $y=-x^k$ for some positive $k,$ we obtain
\begin{equation*}\begin{array}{ll}
\frac{K+1}{K}\left(x^{2k}+3x^4+4x^{3+k}\right)&\le 4x^4+4x^{3+k}+2x^{k+1}+2x^{2k}\\
&\le (K+1)\left(x^{2k}+3x^4+4x^{3+k}\right).\\
\end{array}
\end{equation*}
If $k\in ]1,3[,$ and $x\searrow 0,$ then the previous inequalities fail.}
\end{example}

We are now in the position to show that a function $\varphi\in E^{\mathrm{soft}}(n,K)$ enjoys some regularity properties. Set
$$
E^{\mathrm{soft}}(n):=\bigcup_{K>1} E^{\mathrm{soft}}(n,K).
$$
The following result holds:
\begin{theorem}\label{teorema}
Let $\varphi:\R^n\to \R$ be a convex function in $E^{\mathrm{soft}}(n).$ Then
\begin{itemize}

\item[i.] $\varphi\in {\mathcal C}^1(\R^n)$;

\item[ii.] if, in addition, $\varphi$ has bounded sections, then it is strictly convex.
\end{itemize}
\end{theorem}
{\it Proof}\quad i. By the well known characterization of differentiability for a real-valued convex function (see, for instance, Theorem 25.1 in [7]), $\varphi$ is differentiable in $\R^n$ if and only if, for every $x\in \R^n,$ $\partial \varphi(x)$ is a singleton. Suppose, by contradiction, that there exists $\x\in \R^n$ such that $p,q\in\partial\varphi(\x),$ $p\neq q,$ and set $v=q-p.$ Since $\partial \varphi(\x)$ is a convex set, the segment $p+tv,$ with $t\in [0,1],$ is contained in $\partial \varphi(\x).$ If we consider the function
$$
\psi:\R\to \R, \qquad \psi(s)=\varphi(\x+sv),
$$
one can easily see that $(p+tv)\cdot v\in \partial \psi(0),$ for every $t\in [0,1].$
The function $\psi$ is convex and belongs to $E^{\mathrm{soft}}(1,K),$ being the restriction to a line of a function in $E^{\mathrm{soft}}(n,K).$
In particular, $\psi'$ is defined on a dense subset of $\R,$ according to the Rademacher Theorem. Let $s<0$ be a point in the domain of $\psi'.$
Taking into account the inequalities guaranteed by Proposition \ref{prop:inequalities}, for every $t\in [0,1],$ we have that
\begin{equation*}
\begin{array}{ll}
\frac{K+1}{K}\left(\psi(0)-\psi(s)+\psi'(s)s\right)&\le (\psi'(s)-(p+tv)\cdot v) s\\
&\le (K+1)\left(\psi(0)-\psi(s)+\psi'(s)s\right).\\
\end{array}
\end{equation*}
Dividing by $-s,$ we have
\begin{equation*}
\begin{array}{ll}
\frac{K+1}{K}\left(\frac{\psi(0)-\psi(s)}{-s}-\psi'(s)\right)&\le (p+tv)\cdot v-\psi'(s)\\
&\le (K+1)\left(\frac{\psi(0)-\psi(s)}{-s}-\psi'(s)\right),\\
\end{array}
\end{equation*}
for every $t\in [0,1].$
In particular,
\begin{equation*}
\begin{array}{ll}
\frac{K+1}{K}\left(\frac{\psi(0)-\psi(s)}{-s}-\psi'(s)\right)
&\le  p\cdot v-\psi'(s)< p\cdot v-\psi'(s)+\|v\|^2 \\
&\le (K+1)\left(\frac{\psi(0)-\psi(s)}{-s}-\psi'(s)\right),\\
\end{array}
\end{equation*}
and therefore
$$
\left(\frac{K^2-1}{K}\right) \left(\frac{\psi(0)-\psi(s)}{-s}-\psi'(s)\right)
\ge \|v\|^2.
$$
Setting $\frac{K^2-1}{K}=K_0$, and $\frac{\|v\|^2}{K_0}=T,$ we get that
$$
\frac{\psi(0)-\psi(s)}{-s}-\psi'(s)\ge T;
$$
this is equivalent to say that $0\notin S_{\psi}(s,T),$ for every $s<0$ in the domain of $\psi',$ and for $T$ independent of $s.$
However, this cannot occur; otherwise, by the soft engulfing property enjoyed by $\psi,$ one would have $s\notin S_{\psi}(0,(p+tv)\cdot v,T/K),$ for some $t\in [0,1],$ not possible since $S_{\psi}(0,(p+tv)\cdot v,T/K)$ is open, $0$ is in the interior, and $s$ can be taken close enough to $0.$ This proves that, for every $x\in \R^n,$ $\partial\varphi(x)$ is a singleton, and therefore $\varphi$ is differentiable on $\R^n.$
Finally, taking into account Theorem 25.5 in [7], since $\varphi$ is a proper convex differentiable function, the gradient mapping $\nabla \varphi$ is continuous within $\R^n.$

ii. Firstly, note that the graph of $\varphi$ does not contain any half-line, otherwise $\varphi$ would have at least an unbounded section.
Thus, we will argue by contradiction assuming that the graph of $\varphi$ contains a segment, i.e., there exists $x\in \R^n,$ $v\in \R^n,$ such that $(x+sv,a+sb)\in \mathrm{gph}(\varphi)$ if and only if $s\in [-1,0].$  The problem can be reduced to a one-dimensional problem, by considering the function
$$
\psi:\R\to \R,\qquad \psi(s)=\varphi(x+sv).
$$
Since both the convexity and the engulfing property are invariant with respect to perturbations by affine functions, we can assume, without loss of generality, that
$\psi(s)=0$ if and only if $s\in [-1,0].$
 Note that, due to  convexity and  differentiability, $\partial\psi(s)=\{0\}$ for every $s\in [-1,0]$; moreover,
 $$
 \lim_{s\to -\infty}\psi(s)=+\infty,\qquad \lim_{s\to +\infty}\psi(s)=+\infty,
 $$
and, in particular, $\psi$ establishes a one-to-one correspondence of $[0,+\infty)$ with itself.
Let us consider the section $S_\psi(-1,t);$ in particular,
$$
\mathrm{cl}\,{S_{\psi}(-1,t)}\supset \{s>0:\, \psi(s)\le t\}.
$$
Let $s_t$ be the unique positive point in $\mathrm{cl}\,{S_{\psi}(-1,t)}$ such that $\psi(s_t)=t.$ From the soft engulfing property,
$$
-1\in \mathrm{cl}\,{S_{\psi}(s_t,Kt)},
$$
i.e.,
$$
\psi(-1)-\psi(s_t)-\psi'(s_t)(-1-s_t)\le Kt,
$$
or, equivalently,
\begin{equation}\label{eq:strict1}
\psi'(s_t)(1+s_t)\le (1+K)\psi(s_t).
\end{equation}
From the first inequality in (\ref{caratterizzazione con p}), taking $x=0$ and $y=s_t$, we have that
\begin{equation}\label{eq:strict2}
\left(1+\frac{1}{K}\right)\psi(s_t)\le \psi'(s_t)s_t.
\end{equation}
From \eqref{eq:strict1} and \eqref{eq:strict2}, we get that
$$
\left(1+\frac{1}{K}\right)\psi(s_t)+\psi'(s_t)\le (1+K)\psi(s_t),
$$
or, equivalently,
$$
\psi'(s_t)\le \left(K-\frac{1}{K}\right)\psi(s_t).
$$
Since $\psi:[0,+\infty)\to [0,+\infty)$ is one-to-one, the previous inequality holds
for all $s=s_t>0.$
Then,
$$
\frac{\psi'(s)}{\psi(s)}-\left(K-\frac{1}{K}\right)\le 0,
$$
i.e.,
$$
\frac{{d}}{{ds}}\left(\ln\psi(s)-\left(K-\frac{1}{K}\right)s\right)\le 0,\qquad \forall s>0.
$$
This implies that $\ln\psi(s)-\left(K-\frac{1}{K}\right)s$ is decreasing on $(0,+\infty),$ and therefore there exists
$$
\lim_{s\searrow 0} \left(\ln\psi(s)-\left(K-\frac{1}{K}\right)s\right)=L>-\infty.
$$
Since $\lim_{s\searrow 0}\psi(s)=0^+,$ we get a contradiction.
\qed
\begin{remark}
From the theorem above, if $\varphi\in E^{\mathrm{soft}}(n,K)$, then $\partial \varphi(x)=\{\nabla \varphi (x)\}$, for every $x\in\R^n$, and relation \eqref{caratterizzazione con p} becomes, for every $x, y\in\R^n$,
\begin{equation}\label{caratterizzazione con gradiente}
\begin{array}{ll}
\displaystyle\frac{K+1}{K}&\left(\varphi(y)-\varphi(x)-\nabla \varphi (x)\cdot(y-x)\right)\le (\nabla \varphi (x)-\nabla \varphi (y))\cdot(x-y)\\
&\qquad\qquad\qquad\qquad\le (K+1)\left(\varphi(y)-\varphi(x)-\nabla \varphi (x)\cdot(y-x)\right),\\
\end{array}
\end{equation}
exactly as in Theorem 4 in [6].
\end{remark}

\begin{proposition}\label{proposition soft e non soft}
Let $\varphi:\R^n\to \R$ be a convex function with bounded sections. If $\varphi\in E^{\mathrm{soft}}(n,K),$ then $\varphi\in E(n,K'),$ where $K'=2K(K+1)$.
\end{proposition}
{\it Proof}\quad Let $\varphi\in E^{\mathrm{soft}}(n,K)$. From  Theorem \ref{teorema}, this function is differentiable and strictly convex. Thus, by applying Theorem 4 in [6], the assertion follows.
\qed

\section{Equivalence between engulfing and soft engulfing}\label{uguaglianza}
In this section, we will show that the soft engulfing property, apparently less demanding than the engulfing one, requires in fact the same condition on the functions. We will first prove this result in the one-dimensional case.

\begin{proposition}\label{E=E soft 1}
The following equality holds: $
 E^{\mathrm{soft}}(1)=E(1).
$
Moreover, the unique functions defined on $\R,$ with some unbounded sections and enjoying the engulfing property, are affine.
\end{proposition}
{\it Proof}\quad Let $\varphi\in E^{\mathrm{soft}}(1,K);$ we have only to prove that $\varphi\in E(1,K'),$ for some $K'>1$. From Theorem \ref{teorema}-i, $\varphi$ is differentiable everywhere.
Firstly, let us assume that $\varphi$ is an affine function. It is easy to see that $S_\varphi(x,t)=\R,$
for every $x\in\R$ and $t>0;$ thus,
 $\varphi\in E(1,K)$, for every $K>1$.

Let us now assume that $\varphi$ is not affine; we will show that all the sections are bounded and, by Proposition \ref{proposition soft e non soft}, the assertion will follow.
Suppose, by contradiction, that there exists $x_0\in\R,$
and $t>0,$ such that $S_\varphi(x_0,t)$ is unbounded; since the sections are open and convex subsets of $\R$, and since $\varphi$ is not affine, there exists $\overline x \in\R$ such that, either
$$(\overline x,\infty)=S_\varphi( x_0,t),\qquad \texttt{\rm with}\ \overline x>x_0,
$$
or
$$
(-\infty,\overline x)=S_\varphi( x_0,t),\qquad \texttt{\rm with}\ \overline x>x_0.
$$
Let us consider only the first case, since the second one can be treated similarly.
 For every $ x>  x_0$ we have that $x\in S_\varphi( x_0,t)$ and, by the soft engulfing property,  $x_0 \in S_\varphi( x ,Kt).$ Then,
\begin{eqnarray*}
\varphi(x)-\varphi(x_0)-\varphi'(x_0)(x-x_0) &<&t,\\
\varphi(x_0)-\varphi(x)- \varphi'(x)(x_0-x) &<&Kt,
\end{eqnarray*}
for every $x> x_0.$ The two inequalities imply
\begin{equation}\label{pro 0}
(\varphi'(x)- \varphi'(x_0))(x-x_0)<(K+1)t,\qquad \forall \ x> x_0.
\end{equation}
The convexity of $\varphi$ implies that $\varphi'(x)\ge \varphi'(x_0)$.  Taking $x\to+\infty$ in \eqref{pro 0}, we must have $\varphi'(x)=\varphi'(x_0),$ for every $x>x_0$. Therefore, for every $x\ge x_0,$  the equality
$\varphi(x)=\varphi(x_0)+\varphi'(x_0)(x-x_0)$ holds, and $(\overline x,\infty)=S_\varphi( x,t).$
Now, for every $\epsilon>0,$ we have that $\overline x +\epsilon\in S_\varphi( x,t),$ for every $x>x_0$; the soft engulfing property implies that $x \in S_\varphi(\overline x +\epsilon ,Kt).$ Therefore,
\begin{eqnarray*}
\varphi(\overline x +\epsilon)-\varphi(x)-\varphi'(x_0)(\overline x +\epsilon-x)&<&t,\\
\varphi(x)-\varphi(\overline x +\epsilon)-\varphi'(\overline x+\epsilon)(x-\overline x -\epsilon)&<&Kt,
\end{eqnarray*}
for every $\epsilon>0$ and $x\ge x_0$. The two inequalities imply
\begin{equation}\label{pro 1}
(\varphi'(x_0)-\varphi'(\overline x+\epsilon))(x-\overline x -\epsilon)<(K+1)t,\qquad \forall \epsilon>0,\ x\ge x_0.
\end{equation}
The convexity of $\varphi$ gives that $\varphi'(x_0)\ge  \varphi'(\overline x+\epsilon)$. Taking $x\to+\infty$ in \eqref{pro 1}, we obtain that $\varphi'(x_0)= \varphi'(\overline x+\epsilon),$ for every $\epsilon$, i.e. $\varphi'(\overline x+\epsilon)=\varphi'(x_0)$. Hence, since $\varphi$ is continuous,
$\varphi(x)=\varphi(x_0)+\varphi'(x_0)(x-x_0)$ for every $x\ge \overline x$. This implies that $\overline x\in S_\varphi(x_0,t),$ which contradicts the assumption on $S_\varphi(x_0,t)$.
\qed

The proof of the previous result in the case $n>1$ is  more delicate, since there exist convex functions, that are not affine, and whose graph contains lines:
\begin{example}\quad {\rm
Let $\varphi_0:\R\to\R$ be defined by $\varphi_0(x)=x^2$; $\varphi_0\in  E(1,K),$ for some $K$. If we consider the function $\varphi:\R^2\to\R$ defined by $\varphi(x,y)=x^2$, then, clearly, $S_\varphi((x,y),t)=S_{\varphi_0}(x,t)\times \R$. It is easy to verify that $\varphi\in  E(2,K)$.}
\end{example}

Taking into account the result of the previous proposition, we are now in the position to state and prove our second main result, in the case $n\ge 2.$
\begin{theorem}\label{th:eng=eng}
For every positive integer $n,$ the engulfing and the soft engulfing are equivalent properties, i.e., $
 E^{\mathrm{soft}}(n)=E(n).$
\end{theorem}
\psn
{\it Proof}\quad
We will show that, if $\varphi\in  E^{\mathrm{soft}}(n,K)$, then  $\varphi\in E(n,K'),$ for some $K'>1$.
Let $y\in S_{\varphi}(x,t);$ we will show that $S_{\varphi}(y,t)\subseteq S_{\varphi}(x,K't).$

Suppose, first, that $\varphi(0)=0$ and $\nabla \varphi(0)=0,$ and take any $z\in S_\varphi(y,t),$ $z\neq x.$ Let us consider the set
$$
\tilde{I}=\{w\in \R^n:\, w=(1-s)x+sz,\,s\in \R\}\cap S_{\varphi}(y,t),
$$
and denote by $I$ the convex subset of $\R$ defined as
$$
I=\{s\in \R:\, (1-s)x+sz\in\tilde{I}\}.
$$
$I$ can be bounded, or unbounded. From Proposition \ref{E=E soft 1}, if $I$ is unbounded, then $I=\R,$ and the function $s\mapsto \varphi_{x,z}(s)=\varphi((1-s)x+sz)$ is affine, i.e., there exists $q\in \R$ such that $\varphi_{x,z}(s)=\varphi_{x,z}(0)+qs.$ This implies that, for every $h\in \R,$ for every $l>0,$
$$
S_{\varphi_{x,z}}(h,l)=\R.
$$
In particular, if $l=t,$
this implies that $0\in S_{\varphi_{x,z}}(h,t),$ for every $h\in \R.$
From the soft engulfing property, by choosing $h=1,$ one has that $1\in S_{\varphi_{x,z}}(0,Kt),$ i.e.,
$$
 \varphi_{x,z}(1)<\varphi_{x,z}(0)+\varphi'_{x,z}(0)+Kt,
 $$
 or, equivalently,
 $$
 \varphi(z)<\varphi(x)+\nabla\varphi(x)\cdot(z-x)+Kt.
 $$
 This proves that $z\in S_{\varphi}(x,Kt).$

 Suppose now that $I$ is bounded. In this case, by the convexity of $\varphi_{x,z},$ let $s_1<s_2$ be the real numbers such that
 $$
 t=\varphi_{x,z}(s_1)=\varphi_{x,z}(s_2),\qquad I=]s_1,s_2[.
 $$
 By the Rolle Theorem, there exists $s_0\in ]s_1,s_2[$ such that $\varphi'_{x,z}(s_0)=0.$
 Then, for some $l>0,$ one has that
 $$
 I=S_{\varphi_{x,z}}(s_0,l),
 $$
 i.e.,
 $$
 \varphi_{x,z}(s_1)-\varphi_{x,z}(s_0)=l=\varphi_{x,z}(s_2)-\varphi_{x,z}(s_0).
 $$
 The function $\varphi_{x,z}$ is nonnegative, and therefore $l\le t.$
 Since $\varphi_{x,z}(0)=\varphi (x)$, and $x\in S_{\varphi}(y,t),$ we have that
 $$
 0\in S_{\varphi_{x,z}}(s_0,t).
 $$
 Since $\varphi_{x,z}\in E^{\mathrm{soft}}(1,K),$ from Proposition \ref{E=E soft 1} we get that $\varphi_{x,z}\in E(1,K')$, with $K'=2K(K + 1)$,
and $$
 S_{\varphi_{x,z}}(s_0, t)\subseteq S_{\varphi_{x,z}}(0,K't).
 $$
 In particular,  $z\in S_\varphi(y,t)$ implies that $1\in  S_{\varphi_{x,z}}(s_0, t).$
 Thus,
 $$
 1\in S_{\varphi_{x,z}}(0,K't),
 $$
 i.e.,
 $$
 z\in S_\varphi(x,K't).
 $$
 In the general case, we can consider the function $\psi(x)=\varphi(x)-\varphi(0)-\nabla \varphi(0)\cdot x,$ and note that $\psi\in E(n,K)$ $(E^{\mathrm{soft}}(n,K)),$ if and only if $\varphi$ satisfies the same engulfing condition.

 \qed
\section{Conclusions}
The engulfing property of a convex function is a rather technical condition arising in connection with the Monge-Amp\`{e}re measure.
We introduce a relaxed version of that property, called soft engulfing, and prove that it is actually equivalent to the original
property. Additionally, we show that a convex function with the engulfing property is continuously differentiable and, in case it has bounded sections, it is
strictly convex. All this is done by means of purely convex analytic elementary
techniques.
Following this line of investigation, our aim is to study the engulfing property in the subriemannian setting of the Heisenberg group $\mathbb{H}^n,$ by considering the so called $H$-convex functions and the notion of sections for such functions, which takes into account the peculiar geometry of $\mathbb{H}^n $ (see, for instance,  [8]). In this framework, some important results and tools of the classical convex analysis can be fruitfully extended (see, for instance, [9]), and applied to the investigation of the engulfing property in $\mathbb{H}^n$ (see [10]).

\end{document}